\documentclass[11pt,twoside]{article}

\usepackage{amsmath}
\usepackage{amssymb}
\usepackage{fancyhdr}
\usepackage{graphicx}
\usepackage{subfigure}
\usepackage[textwidth=16cm,textheight=22cm,top=3cm,bottom=3cm]{geometry}

\newcommand{\drop}[1]{}

\newcommand{\Div}{\operatorname{div}}

 \newcommand{\qtext}[1]{\quad\text{#1}}
\newcommand{\p}{\partial}
\newcommand{\R}{\mathbb{R}}
\newcommand{\eps}{\varepsilon}
\newcommand{\vfi}{\varphi}
\renewcommand{\O}{\Omega}

\newcommand{\fer}[1]{(\ref{#1})}

\newcommand{\no}{\noindent}
\newcommand{\ud}{u^{\delta}}
\newcommand{\rd}{r^{\delta}}
\newcommand{\Pd}{{\text{(P)}$_\delta$}}
\newcommand{\Ps}{{\text{(PS)}$_\delta$}}

\newcommand{\grad}{\nabla}
\newcommand{\abs}[1]{| #1 |}
\newcommand{\nor}[1]{\| #1 \|}

\newcommand{\wto}{\rightharpoonup}

\newtheorem{theorem}{Theorem}
\newtheorem{remark}{Remark}

\newcommand{\q}{\mathbf{q}}

\title{Analysis of a splitting-differentiation population model leading to cross-diffusion
\thanks{Supported by Spanish MCI Project MTM2013-43671-P.}}

\author{Gonzalo Galiano  \thanks{Dpt. of Mathematics, Universidad de Oviedo,
 c/ Calvo Sotelo, 33007-Oviedo, Spain ({\tt galiano@uniovi.es, selgasvirginia@uniovi.es})}
    \and Virginia Selgas\footnotemark[2] }
\date{}

\begin{document}
\maketitle

\begin{abstract}
Starting from the dynamical system model capturing the splitting-differentiation process of populations, we extend this notion to show how the speciation mechanism from a single species leads to the consideration of several well known evolution cross-diffusion partial differential equations.

Among the different alternatives for the diffusion terms, we study the model introduced by Busenberg and Travis, for which we prove the existence of solutions in the one-dimensional spatial case. Using a direct 
parabolic regularization technique, we show that the problem is well posed in the space of bounded variation functions, and demonstrate with a simple example that this is the best regularity expected for solutions.

We numerically compare our approach to other alternative regularizations previously introduced in the literature, for the particular case of the contact inhibition problem. Simulation experiments indicate that the numerical scheme arising from the approximation introduced in this article outperforms those of the existent models from the stability point of view.

\vspace{0.25cm}

\no\emph{Keywords: }Cross-diffusion system, population dynamics, contact inhibition, splitting and differentiation,
numerical simulations, existence of solutions.
\end{abstract}

% \MSC 35K55 \sep 65M30 \sep 92D25

 \section{Mathematical model and main result}

 \subsection{The splitting-differentiation model in terms of ODEs}
 In \cite{sanchez-palencia}, S\'anchez-Palencia analyzes the situation in which a species,  with population density $U$,  splits, due to a number of factors, into two different species with population densities $U_1$ and $U_2$, but still keeping the original ecological behavior. 
 
 Thus, we assume that $U$ followed a logistics law until time $t^*$, i.e.
 \begin{equation}
 \label{sp.1}
  U'(t)= U(t) \big(\alpha -\beta U(t)\big),\qtext{for } t\in(0,t^*), \quad U(0)=U_0>0.
 \end{equation}
Then, after splitting, $(U_1,U_2)$ is assumed to satisfy the Lotka-Volterra system
\begin{align}
\label{sp.2}
  U_i'(t)= U_i(t) \big(\alpha -\beta \big(U_1(t)+U_2(t)\big)\big),\qtext{for } t\in(t^*,T), \quad U_i(t^*)=U_{i0}>0,
\end{align}
for $i=1,2$,
with $U_{10}+U_{20}=U(t^*)$. Note that, under this splitting \emph{without differentiation}, $U_1+U_2$ still satisfies \fer{sp.1}
for $t\geq t^*$.
 
Although problem \fer{sp.1} has a unique non-trivial equilibrium, $U_{\infty}=\alpha /\beta$, problem \fer{sp.2} has a continuum set of non-trivial equilibria given by
all the combinations of $U_{1\infty}\geq0$ and $U_{2\infty}\geq0$  such that $U_{1\infty}+U_{2\infty}=\alpha /\beta$. 

In \cite{sanchez-palencia}, the author analyzes how the differentiation of populations $U_1$, $U_2$ after  splitting, understood as a perturbation in the Lotka-Volterra coefficients, affects  the equilibrium of the system. By splitting \emph{with differentiation} we mean that $(U_1,U_2)$ is a solution of the following problem:
\begin{align}
\label{sp.3}
  U_i'(t)= U_i(t) \big(\alpha_i - \big(\beta_{i1}U_1(t)+\beta_{i2}U_2(t)\big)\big),\qtext{for } t\in(t^*,T), \quad U_i(t^*)=U_{i0}>0,
\end{align}
for $i=1,2$, with $U_{10}+U_{20}=U(t^*)$. Observe that, under this splitting, $U_1+U_2$ does not satisfy, in general, problem \fer{sp.1} for $t\geq t^*$.

The main conclusion of \cite{sanchez-palencia} is that the differentiation mechanism selects, in general, a unique solution for the equilibrium system, having therefore a stabilizing effect.
 
\subsection{The splitting model in terms of PDEs: cross-diffusion}
In this article, we extend the previous dynamical system models to the case of space dependent  population densities.
We start considering the dynamics of one single species population
satisfying
\begin{equation}
\left\{
\begin{array}{ll}
\partial_t u -\Div J(u) =f(u) & \qtext{in } Q_{(0,t^*)},\\
J(u)\cdot \boldsymbol{\nu}=0 & \qtext{on } \Gamma_{(0,t^*)}, \\
u(0,\cdot)=u_0 \geq 0 & \qtext{on }\O,
\end{array}
\right.
\label{p.u}
\end{equation}
where $\O\subset \R^N$ is a bounded domain with Lipschitz continuous boundary, $\p\O$,  $Q_{(0,t^*)}=(0,t^*)\times\O$, and $\Gamma_{(0,t^*)}=(0,t^*)\times\p\O$ is the parabolic boundary of $Q_{(0,t^*)}$. The vector $\boldsymbol{\nu}$ is the outwards canonical normal to $\p\O$.
The growth-competition term is assumed to have the logistic form 
$f(u)=u(\alpha-\beta u)$, and the flow to be given by
$ J(u)=u\grad u+u\q$.

As it is well known, the term
$u\grad u$ captures the individuals aversion to
overcrowding, while $\q$ is usually determined by an environmental potential, $\q=-\grad \vfi$, whose minima represent attracting points for the populations. For biological background and origins of the model see, for instance, \cite{okubo}.

After splitting, the new two populations, $u_1$ and $u_2$, satisfy, for $i=1,2$, 
\begin{equation}
\left\{
\begin{array}{ll}
\partial_t u_i -\Div J_i(u_1,u_2) =f_i(u_1,u_2) & \qtext{in } Q_{(t^*,T)},\\
J_i(u_1,u_2)\cdot \boldsymbol{\nu}=0 & \qtext{on } \Gamma_{ (t^*,T)}, \\
u_i(t^*,\cdot)=u_{i0} & \qtext{on }\O,
\end{array}
\right.
\label{p.u2}
\end{equation}
 with $u_{i0}$ such that $u_{10}+u_{20}=u(t^*,\cdot)$, and with $J_i$ and $f_i$ to be defined. 

Assuming, like in the dynamical system model, that the possible differentiation process only takes place through the growth and the inter- and intra-competitive behavior of the new species implies that the split flows must satisfy $J_1(u_1,u_2)+J_2(u_1,u_2)=J(u_1+u_2)$, that is
\[
 J_1(u_1,u_2)+J_2(u_1,u_2)=(u_1+u_2)\grad(u_1+u_2) + (u_1+u_2)\q.
\]
Being clear the  way of defining the linear transport term of $J_i$, the nonlinear diffusive term 
admits several reasonable decompositions. For instance, in \cite{Galiano2012}, the following  splitting was considered
\begin{equation}
 \label{def.flows2}
J_i(u_1,u_2)=u_i\grad u_i + b_i \grad(u_1u_2)+u_i\q ,
\end{equation}
with $b_i\geq 0$, and $b_1+b_2 = 1$. Under this splitting, 
problem \fer{p.u2} takes the form of the cross-diffusion model introduced by Shigesada et al. \cite{SKT79}, for which a thorough mathematical analysis does exist, see for instance  \cite{ggj2,gambino,andreianov11,berres11,Ruiz-Baier2012,Juengel2014} for numerical approaches, \cite{chen,gv,Gambino2012,Gambino2013,Desvillettes2014} for analytical and qualitative results, or \cite{gjv,Galiano2011,Juengel2012} for applications.

In this paper, we consider the alternative splitting
\begin{equation}
 \label{def.flows}
J_i(u_1,u_2)=u_i\grad (u_1+u_2)+ u_i\q, 
\end{equation}
which brings problem \fer{p.u2} to the form of the  Busenberg and Travis model \cite{busenberg83}. 
Although apparently simpler than \fer{def.flows2}, no general proof of existence of solutions does exist for problem \fer{p.u2} with flows given by \fer{def.flows}.
However, some partial results related to the cell-growth contact-inhibition problem may be found in \cite{bertsch85,bertsch2010,bertsch2012,gsv}, as well as in \cite{Galiano2014,gs} for other
specific situations.

\subsection{Differentiation after splitting}
 According to the species behavior after splitting, we consider two problems arising from two different sets of Lotka-Volterra terms:
\begin{align}
& f_i(u_1,u_2)=u_i(\alpha-\beta (u_1+ u_2)),& &\text{(non-differentiation)}& \label{f.nm}\\
 &f_i(u_1,u_2)=u_i \big(\alpha_i - \big(\beta_{i1}u_1+\beta_{i2}u_2\big)\big),& &\text{(differentiation)}\label{f.m}&
\end{align}
for $i=1,2$. We shall refer to problem \fer{p.u2} with flows given by \fer{def.flows} and with $f_i$ given by \fer{f.nm} and \fer{f.m} as  to \textbf{problems (ND) and (D)}, respectively. Observe that these are the PDE versions that generalize the non-differentiation and differentiation ODE problems \fer{sp.2} and \fer{sp.3}, introduced in \cite{sanchez-palencia}.

The existence of solutions of problem (ND) was proven in \cite{Galiano2014} in the multi-dimensional case. The proof is based on the construction of solutions, $(\ud_1 ,\ud_2)$, to a nonlinear parabolic regularization of problem (ND)  such that, although $\ud_i$ only converges weakly to the solution 
of the limit problem, $u_i$, the global density,  $\ud_1+\ud_2$, converges strongly to $u_1+u_2$. Thus, weak and strong convergences are compensated so that the limit problem (ND) is shown to have a solution. 

However, since no a.e. convergence of $\ud_i$ to $u_i$ was proven in \cite{Galiano2014}, the case of differentiated Lotka-Volterra terms, i.e. problem (D), may not be handled directly with this technique. 

Previously to \cite{Galiano2014},
in the one-dimensional setting and for $q=0$, Bertsch et al. \cite{bertsch2010} proved the existence of a solution  of  problem (D) with the form 
\begin{equation*}
%\label{sol.b}
  u_1(t,x)=r(t,x)u(t,x),\quad u_2(t,x)=\big(1-r(t,x)\big) u(t,x),
\end{equation*}
where $u=u_1+u_2$ and $0\leq r \leq 1$ solve certain parabolic-hyperbolic auxiliary problem, see  problem (P)$_B$ in Section~\ref{sec:app}.

Their proof is based on the parabolic regularization of the auxiliary problem  %of \fer{p.u4}
together with the use of the Lagrangian flows (characteristics) associated to $\p_x u$. In particular,
 and important for the results proven in this article, 
they obtained \emph{strong convergence} of the sequence of approximated solutions
to a limit, which is identified as a solution of  problem~(D).

The first aim of this  paper is to show the existence of solutions of problem (D) with an alternative proof to that given in  \cite{bertsch2010}. Our proof is  based on the direct parabolic regularization used for problem (ND) in \cite{Galiano2014}, and  takes advantage of the techniques  employed in \cite{bertsch2010} to obtain the strong convergence of each component of the regularized problem. 

In fact, like in \cite{bertsch2010}, the space regularity of solutions of problem (D) is shown to be of bounded variation, $BV(\O)$, which seems to be the optimal expected regularity, see Theorem~\ref{th.easy} in Section~\ref{sec:numerics}, for an example.  

More explicitly, we prove the following result in Section~\ref{sec:proof}. Here, we retake the usual notation $[0,T]$ for the time domain, replacing $[t^*,T]$. 
 \begin{theorem}
\label{th.1}
 Let $\O\subset\R$ be a bounded interval and $T>0$ be arbitrarily fixed. Let $q\in C^{0,1}(\bar Q_T)$, with $q(t,\cdot)=0$ on $\p\O$ for all $t\in[0,T]$,  and $f_i:\R^2\to\R$ be given by \fer{f.m}. Assume that $u_{10},u_{20}\in BV(\O)$ are non-negative, with $u_0:=u_{10}+u_{20}>0$ in $\bar\O$,   $u_0\in C^{0,1}(\bar \O)$ and $\p_x u_0=0$ on $\p\O$.  Then, there exists a weak solution of problem (D), $(u_1,u_2)$, with, for $i=1,2$, 
 \begin{enumerate}
 \item $u_i \geq0$ a.e. in $Q_T$.
  \item $u_i \in  L^\infty(0,T;BV(\O))\cap BV(0,T;L^1(\O))$.
  \item $u:=u_1+u_2 \in  C^{0,1}(\bar Q_T) \cap L^2(0,T;H^2(\O))$.
  \item For all $\vfi \in C^1(0,T;H^1(\O))$ with $\vfi(T,\cdot)=0$ in $\O$, 
  \begin{align*}
 \int_{Q_T} u_i\p_t \vfi = \int_{Q_T} \big(u_i\p_x u + u_iq \big) \p_x\vfi- \int_{Q_T} f_i(u_1,u_2)\vfi- \int_\O u_{i0}\vfi(0,\cdot).
 %\label{def.wd}
\end{align*}
\end{enumerate}
\end{theorem}

The second aim of this article is to numerically investigate and compare the resulting Finite Element schemes of each regularizing approach, see  problems \Pd~ and (P)$_B$ in Section~\ref{sec:app}.  
We focus our attention into two model problems: (i) the \emph{invasion} problem, arising in Tumor theory, as suggested by Bertsch et al. \cite{bertsch2010}, in which an initial small perturbation (tumor) inside a healthy tissue evolves in time keeping a sharp interface; and (ii) a model problem with a Barenblatt-based explicit solution, for which a detailed comparison to the approximated solutions is given, see Section~\ref{sec:numerics}.

After the proofs of our results in Section~\ref{sec:proof}, we conclude the article with some 
conclusions, in Section~\ref{sec:conclusions}.

\section{Approximated problems and discretization}\label{sec:app}

In this section we compare the approximated solutions to problem (D) constructed via the regularization scheme employed in the proof of Theorem~\ref{th.1}, and those corresponding to the  approximated problem introduced in \cite{bertsch2010}. The formulation of these problems is the following, respectively. 
\begin{equation*}
\text{(P)$_\delta$}\left\{
\begin{array}{ll}
\p_t u_i  - \p_x \big(u_i\p_x (u_1+u_2) + u_iq \big) - \frac{\delta}{2}\p_{xx}(u_i (u_1+u_2)) 
 = f_i(u_1,u_2) & \text{in } Q_{T},\\
\p_x u_i=0 & \text{on } \Gamma_{T}, \\
u_i(0,\cdot)=\ud_{i0} & \text{in }\O,
\end{array}
\right.
%\label{p.nuestro}
\end{equation*}
for $\delta>0$, and some non-negative $\ud_{i0}\in C^1(\bar\O)$ such that $\ud_{i0}\to u_{i0}$ strongly in $BV(\O)$, for $i=1,2$.

\begin{equation*}
\text{(P)$_B$}\left\{
\begin{array}{ll}
\partial_t u -\p_x(u (\p_x u +q)) =F_1(u,r) & \text{in } Q_{T},\\
\partial_t r -(\p_x u+q) \p_x r -\delta_B \p_{xx} r=F_2(u,r) & \text{in } Q_{T},\\
\p_x u =\p_x r=0 & \text{on } \Gamma_{T}, \\
u(0,\cdot)=u_{0},\quad r(0,\cdot)=r_0^{\delta_B} & \text{in }\O,
\end{array}
\right.
%\label{p.u4}
\end{equation*}
for $\delta_B>0$, and for some $r_0^{\delta_B}\in C^1(\bar\O)$ such that $r_0^{\delta_B}\to r_{0}$ strongly in $BV(\O)$, and with 
 \begin{align*}
 & F_1(u,r)=f_1(ru,(1-r)u)+f_2(ru,(1-r)u), \\
 & F_2(u,r)=r(1-r)\Big(\frac{f_1(ru,(1-r)u)}{ru}-\frac{f_2(ru,(1-r)u)}{(1-r)u}\Big).
 \end{align*}
 Recall that, according to \cite{bertsch2010}, for any $\delta_B>0$ there exists a regular solution to problem (P)$_B$, $(u^{\delta_B},r^{\delta_B})$, and that this sequence  converges strongly in $L^1(Q_T)$, as $\delta_B\to0$, to some $(u,r)$ such that   $u_1=ru$, and $u_2=(1-r)u$ are a solution of problem (D).

For the numerical discretization of problems (P)$_\delta$ and (P)$_B$, we consider a fully discrete approximation using finite elements in space and backward finite differences in time. 
The proof of the convergence of the numerical scheme for problem (P)$_\delta$ may be found in \cite{Galiano2014}.

We consider  a quasi-uniform mesh on the interval $\Omega$, $\{\mathcal{T}_h\} _h$, with $h$ representing step size.
We introduce the finite element space of continuous  $\mathbb{P}_1$-piecewise elements:
$$
S^h = \{ \chi\in \mathcal{C}(\overline{\Omega} ) ; \, \chi |_\kappa \in\mathbb{P}_1\,\text{ for all } \kappa \in\mathcal{T}_h \} .
$$
The Lagrange interpolation operator is denoted by $\pi ^h : \mathcal{C}(\overline{\Omega} ) \to S^h$. We also introduce 
the discrete semi-inner product on $\mathcal{C}(\overline{\Omega} ) $ and its induced discrete seminorm:
$$
(\eta_1,\eta_2)^h= \int_{\Omega} \pi^h(\eta_1\eta_2) ,\quad |\eta|_h=\sqrt{(\eta,\eta)^h}.
$$
For each $\eps\in (0,1)$ we consider the  function 
\begin{equation}
 \label{landa}
 \lambda_\eps(s)=\left\{
 \begin{array}{ll}
  \eps & \text{if }s\leq \eps,\\
  s & \text{if } \eps\leq s\leq \eps^{-1},\\
  \eps^{-1} & \text{if } s\geq \eps^{-1},
 \end{array}
\right.
\end{equation}
and the linear operator $\Lambda _{\eps} : S^h\to L^{\infty}(\Omega)$  given by, for 
$x_{mp}=(x_j+x_{j+1})/2$
\begin{equation}
\label{landamay}
\Lambda_\eps(z^h) = \lambda_{\eps}(z^h(x_{mp})),  \qtext{in } (x_j,x_{j+1}) .
\end{equation}
% \begin{equation*}
%  \Lambda_\eps(z^h)|_\kappa=\left\{
%  \begin{array}{ll}
%   \dfrac{1}{\lambda_{\eps}(z^h(\xi))} & \text{if }z^h(x_j)\neq z^h(x_{j+1}),\\
%   \dfrac{1}{\lambda_{\eps}(z^h(x_k))} & \text{if }z^h(x_j)= z^h(x_{j+1}),
%  \end{array}
% \right.
% \end{equation*}
% for some  $\xi\in \kappa$, see \cite{grun00}. In the experiments, we always take $\xi=(x_j+x_{j+1})/2$.

For the time discretization, we take in the experiments a uniform partition of $[0,T]$ of time step $\tau$.
For $t=t_0=0$, set $u_{i\eps }^0=u_i^0$, for $i=1,2$.
Then, for $n\geq 1$ the full discretization of problem (P)$_\delta$ reads: Find $u_{i\eps }^{n}\in S^h$ such that
\begin{align*}
%\label{eq:pde_discr.s4}
\tfrac{1}{\tau}\big( u^n_{i\eps }-u^{n-1}_{i\eps } , \chi )^h
 + (1+\tfrac{\delta}{2})\big(U_{i\eps}^n \p_x ( u^n_{1\eps }+u^n_{2\eps })   ,\p_x \chi \big) 
+ \tfrac{\delta}{2}\big((U_{1\eps}^n+U_{2\eps}^n) \p_x u^n_{i\eps }   ,\p_x \chi \big) \\
%\hspace*{3cm}
 +  \big(  \pi^h(q) U_{i\eps}^n   ,\p_x \chi \big)
%\\ [2ex]\hspace*{1cm} 
= \big( \alpha_i u^n_{i\eps } - \lambda_\eps(u^n_{i\eps })(\beta_{i1}\lambda _{\eps } (u^{n-1}_{1\eps })+\beta_{i2}\lambda _{\eps } (u^{n-1}_{2\eps })), \chi \big)^h , %\qquad\forall \chi\in S^h .
\end{align*}
for every $ \chi\in S^h $, where we introduced the notation $U_{i\eps}^n=\Lambda_\eps(u_{i\eps}^n)$.

Similarly, the full discretization of problem (P)$_B$ reads: Set $(u_{\eps }^0, r_\eps^0)=(u_0, r_0)$. 
Then, for $n\geq 1$, find $(u_{\eps }^{n}, r_\eps^n) \in S^h\times S^h$ such that
\begin{align*}
%\label{eq:pde_discr.s5}
\tfrac{1}{\tau}\big( u^{n}_{\eps }-u^{n-1}_{\eps } , \chi )^h
& + \big(U^{n}_{\eps } \p_x u^{n}_{\eps }   ,\p_x \chi \big) 
+  \big(  \pi^h(q) U^{n}_{\eps }   ,\p_x \chi \big) 
 = \big( F_{1\eps}(u^{n}_{\eps },r^{n}_{\eps }), \chi \big)^h ,\\[2ex]
%%%%%%%%%%%%%%%%%%%%%%%%
\tfrac{1}{\tau}\big( r^{n}_{\eps }-r^{n-1}_{\eps } , \chi )^h
& + \delta_B \big( \p_x r^{n}_{\eps }  ,\p_x \chi \big) 
-  \big(  (\p_x u^{n}_\eps + \pi^h(q)) \p_x r^{n}_\eps  ,\p_x \chi \big) \\
&  = \big( F_{2\eps}(u^{n}_{\eps },r^{n}_{\eps }), \chi \big)^h .
\end{align*}
for every $ \chi\in S^h $, where $U_{\eps}^n=\Lambda_\eps(u_{\eps}^n)$, and 
$F_{i\eps}(s,\sigma)= F_i(\lambda _{\eps } (s),\lambda _{\eps } (\sigma))$.

Since the above systems are nonlinear algebraic problems, we use a fixed point argument to 
approximate their solution at each time slice $t=t_n$, from the previous
approximation at $t=t_{n-1}$. Thus, for problem (P)$_\delta$, let $u_{\eps i}^{n,0}=u_{\eps i}^{n-1}$. 
Then, for $k\geq 1$ the problem is to find $u_{\eps i}^{n,k} \in S^h$ such that for 
$i=1,2$, and for all $\chi \in S^h$ 
\begin{align*}
%\label{eq:pde_discr.s4}
\tfrac{1}{\tau}\big( u^{n,k}_{i\eps }-u^{n-1}_{i\eps } , \chi )^h 
&+ (1+\tfrac{\delta}{2})\big(U_{i\eps}^{n,k-1} \p_x ( u^{n,k}_{1\eps }+u^{n,k}_{2\eps }),\p_x \chi \big)\\
& + \tfrac{\delta}{2}\big((U_{1\eps}^{n,k-1}+U_{2\eps}^{n,k-1}) \p_x u^{n,k}_{i\eps }   ,\p_x \chi \big) 
 +  \big(  \pi^h(q) U_{i\eps}^{n,k-1}   ,\p_x \chi \big) \\
& = \big( \alpha_i u^{n,k}_{i\eps } - \lambda_\eps(u^{n,k-1}_{i\eps })(\beta_{i1}\lambda _{\eps } (u^{n-1}_{1\eps })+\beta_{i2}\lambda _{\eps } (u^{n-1}_{2\eps })), \chi \big)^h .
\end{align*}
We then use the stopping criterion $\max _{i=1,2} \nor{u_{\eps,i}^{n,k}-u_{\eps,i}^{n,k-1}}_\infty <\text{tol}$,
for values of $\text{tol}$ chosen empirically, and set $u_i^n=u_i^{n,k}$.
% In some of the experiments we integrate in time until a numerical stationary solution, 
% $u_i^S$, is achieved. This is determined 
% by $\max _{i=1,2} \nor{u_{\eps,i}^{n,1}-u_{\eps,i}^{n,0}}_\infty <\text{tol}_S$,
% where $\text{tol}_S$ is chosen also empirically.

Similarly, for problem (P)$_B$ we use the scheme 
\begin{align*}
%\label{eq:pde_discr.s5}
\tfrac{1}{\tau}\big( u^{n,k}_{\eps }-u^{n-1}_{\eps } , \chi )^h
& + \big(U^{n,k-1}_{\eps } \p_x u^{n,k}_{\eps }   ,\p_x \chi \big) 
+  \big(  \pi^h(q) U^{n,k-1}_{\eps }   ,\p_x \chi \big) \\
& = \big( F_{1\eps}(u^{n,k-1}_{\eps },r^{n,k-1}_{\eps }), \chi \big)^h ,\\[2ex]
%%%%%%%%%%%%%%%%%%%%%%%%
\tfrac{1}{\tau}\big( r^{n,k}_{\eps }-r^{n-1}_{\eps } , \chi )^h
& + \delta_B \big( \p_x r^{n,k}_{\eps }  ,\p_x \chi \big) 
-  \big(  (\p_x u^{n,k-1}_\eps + \pi^h(q)) \p_x r^{n,k-1}_\eps  ,\p_x \chi \big) \\
&  = \big( F_{2\eps}(u^{n,k-1}_{\eps },r^{n,k-1}_{\eps }), \chi \big)^h ,
\end{align*}
with an analogous stopping criterion than above.

\section{Numerical simulations}%: An invasion and the contact inhibition problem}
\label{sec:numerics}

In this section, we present numerical simulations for two sets of initial and Lotka-Volterra
data aiming to clarify the advantages of approximating the solutions of problem (D) either by
the scheme derived from the regularized parabolic-hyperbolic formulation (P)$_B$ of \cite{bertsch2010}, or by the direct viscosity approximation introduced in this article, (P)$_\delta$.

The general conclusion is that approximating the solution, $(u_1,u_2)$, of problem (D)
by the sequence $(u_{1\delta},u_{2\delta})$ of solutions of problem (P)$_\delta$ 
is more robust against instabilities produced in discontinuity points than that given by the sequence $(u_{\delta_B},r_{\delta_B})$ corresponding to  (P)$_B$. Although we only show the results
for two standard examples, we confirmed this conclusion for a variety of data problem.

Parameter data was fixed as follows.  In all the examples, we consider the spatial domain $\O=(-2,2)$ and investigate the mesh sizes $h=0.04$, $0.013$ and $0.008$, whereas the time step is empirically chosen as $\tau=1.e-3$ (in the first experiment) and $1.e-4$ (in the second). 
The regularization parameters are taken as $\delta=h^2$, $\delta_B=2h^2$. To avoid negative values of the discrete solutions, we use the functions $\lambda_\eps$ and $\Lambda_\eps$, see \fer{landa} and \fer{landamay}, respectively, setting $\eps=1.e-10$. Finally, the tolerance for the fixed point stopping criterion is set to $\text{tol}=1.e-8$; in all the examples, we observed good convergence properties of the fixed point algorithm, reaching the prescribed tolerance in less than ten iterations. %Concerning the stationary state tolerance, it is fixed as $1.e-11$ and never reached.

In all the cases, the solutions corresponding to problem (P)$_B$ produce more oscillations than those of problem (P)$_\delta$. To measure these oscillations, we compute an approximation to the zero-crossings of $\p_x u(t,\cdot)$ in $\O$ by 
\[
\textrm{osc}(u)(t)= h\sum \abs{\Delta(\textrm{sign}(\Delta u(t,\cdot)))},
\]
where $\Delta$ stands for the difference between two consecutive node values, and the sum 
runs over the spatial nodes. In fact, we may observe that while solutions of  (P)$_\delta$ attenuate the oscillations when $t$ increases, those corresponding to solutions of  (P)$_B$ are always above some threshold value. In addition,  in all the experiments we observe  that the size of oscillations diminish according to the mesh size, for both approximations.

%Of particular interest, we consider both an invasion problem and a particular example of the contact inhibition problem where the explicit solution can be built.

\subsection{Experiment 1: Invasion}

\begin{figure}[h] 
\centering 
{\includegraphics[width=3cm,height=3cm]{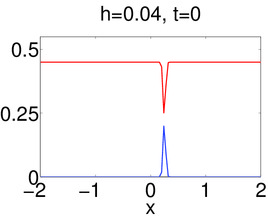}}
{\includegraphics[width=3cm,height=3cm]{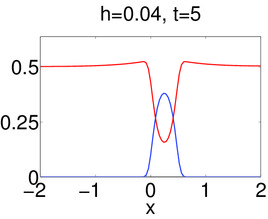}}
{\includegraphics[width=3cm,height=3cm]{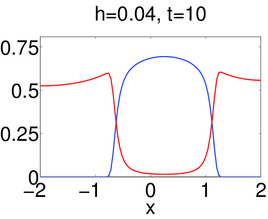}}
{\includegraphics[width=3cm,height=3cm]{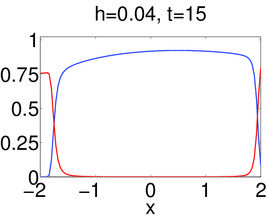}}\\
{\includegraphics[width=3cm,height=3cm]{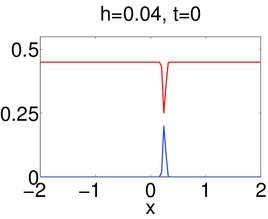}}
{\includegraphics[width=3cm,height=3cm]{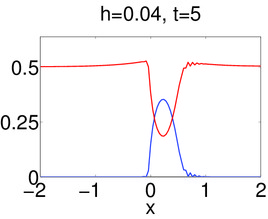}}
{\includegraphics[width=3cm,height=3cm]{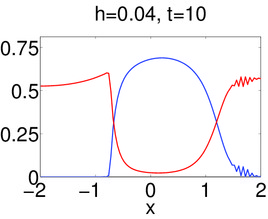}}
{\includegraphics[width=3cm,height=3cm]{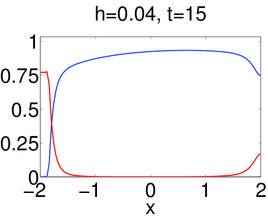}}\\
{\includegraphics[width=3cm,height=3cm]{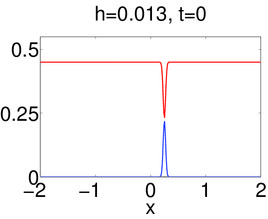}}
{\includegraphics[width=3cm,height=3cm]{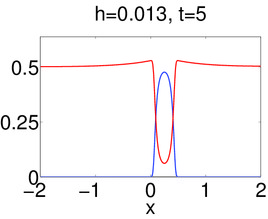}}
{\includegraphics[width=3cm,height=3cm]{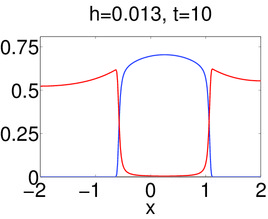}}
{\includegraphics[width=3cm,height=3cm]{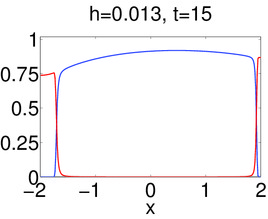}}\\
{\includegraphics[width=3cm,height=3cm]{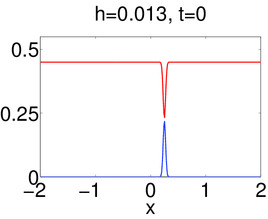}}
{\includegraphics[width=3cm,height=3cm]{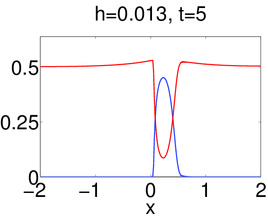}}
{\includegraphics[width=3cm,height=3cm]{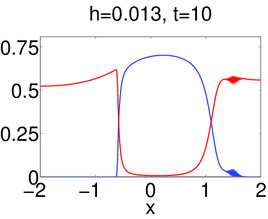}}
{\includegraphics[width=3cm,height=3cm]{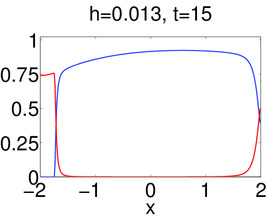}}\\
{\includegraphics[width=3cm,height=3cm]{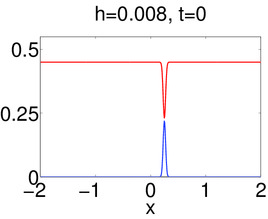}}
{\includegraphics[width=3cm,height=3cm]{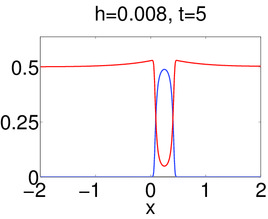}}
{\includegraphics[width=3cm,height=3cm]{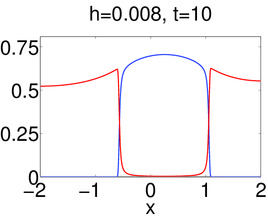}}
{\includegraphics[width=3cm,height=3cm]{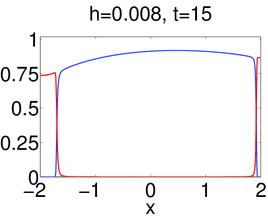}}\\
{\includegraphics[width=3cm,height=3cm]{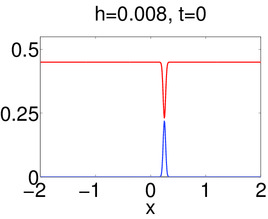}}
{\includegraphics[width=3cm,height=3cm]{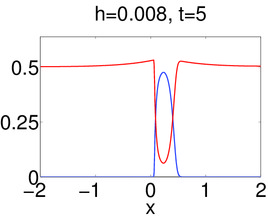}}
{\includegraphics[width=3cm,height=3cm]{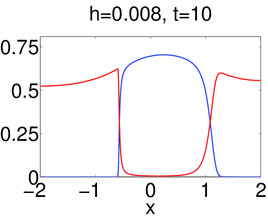}}
{\includegraphics[width=3cm,height=3cm]{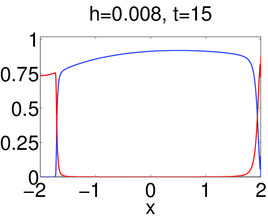}}\\
\caption{  Experiment 1.  
Odd rows correspond to the approximated solution of problem (P)$_\delta$ for several time slices, $t$, while even rows correspond to the solution of problem (P)$_B$.
Rows correspond to different mesh sizes, captured by parameter $h$. Mind the different vertical scales among time slices.}
\label{fig_invasion_sol} 
\end{figure}

In this example we investigate the question (Q2) stated in \cite{bertsch2010}, concerning the invasion of a population (which simulates mutated abnormal cells) over an initially dominant population (which represents the normal cell). More precisely, we take the initial data
$$
u_{10}(x) = 0.22\, \mathrm{exp} (-(x-0.25)^2/0.001)\, , \quad u_{20}(x) = 0.45-u_{10}(x) \qquad\text{for }x\in \O \, .
$$
Besides, we deal with the Lotka-Volterra coefficients
$$
\alpha_i =1 \, , \quad \beta _{ij}= i \qquad\text{for }i,j=1,2 \, .
$$
We conducted the experiments for several choices of the mesh size, $h$. Only the experiments corresponding to $h=0.04$ (101 nodes), $h=0.013$ (301 nodes), and $h=0.008$ (501 nodes) are shown in Figure~\ref{fig_invasion_sol}.

\begin{figure}[h] 
\centering 
{\includegraphics[width=4cm,height=4cm]{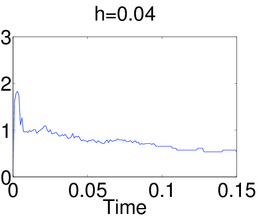}}
{\includegraphics[width=4cm,height=4cm]{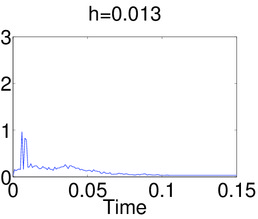}}
{\includegraphics[width=4cm,height=4cm]{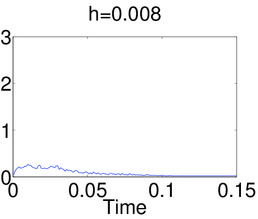}} \\[2ex]
{\includegraphics[width=4cm,height=4cm]{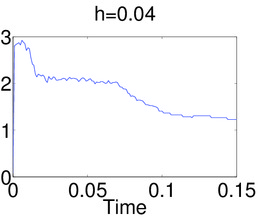}}
{\includegraphics[width=4cm,height=4cm]{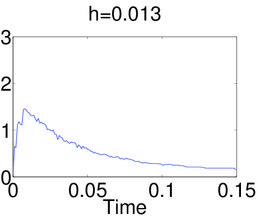}}
{\includegraphics[width=4cm,height=4cm]{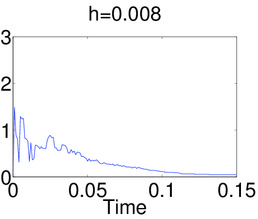}}
\caption{  Experiment 1. Each row corresponds to the oscillation measure $\textrm{osc}(u)(t)$ of the solutions approximated by (P)$_\delta$ and (P)$_B$, respectively. 
Columns correspond to different mesh sizes, captured by parameter $h$.}
\label{fig_invasion_osc} 
\end{figure}

We may check that for the finer mesh (last two rows), both approximations give similar results. 
However, instabilities are already present in the approximation obtained from problem (P)$_B$ for the medium size mesh (row four), which appear amplified for the coarsest mesh (row two). 
For the printed resolution, instabilities are not visible in the approximation 
obtained from problem (P)$_\delta$ (rows one and three). This graphical evidences 
are confirmed through the oscillation measure $\textrm{osc}(u)(t)$, plotted 
in Figure \ref{fig_invasion_osc}.

% \begin{figure}[H] 
% \centering 
% {\includegraphics[width=5cm,height=4cm]{figINV_100_osc.jpg}}
% {\includegraphics[width=5cm,height=4cm]{figINV_UR_100_osc.jpg}}\\[2ex]
% {\includegraphics[width=5cm,height=4cm]{figINV_300_osc.jpg}}
% {\includegraphics[width=5cm,height=4cm]{figINV_UR_300_osc.jpg}}\\[2ex]
% {\includegraphics[width=5cm,height=4cm]{figINV_500_osc.jpg}}
% {\includegraphics[width=5cm,height=4cm]{figINV_UR_500_osc.jpg}}
% \caption{  Experiment 1. Each column corresponds to the oscillation measure of the solutions approximated by (P)$_\delta$ and (P)$_B$, respectively. 
% Rows correspond to different mesh sizes, captured by parameter $h$.}
% \label{fig_invasion_osc} 
% \end{figure} 

\subsection{Experiment 2: A Barenblatt-based explicit solution} 

The contact inhibition problem is a particular case of problem (D) arising in tumor dynamics theory, see \cite{Chaplain2006} for the modeling background. In this problem, the initial density distributions  are segregated, satisfying
\begin{equation}
\label{ci.id} \operatorname{supp} u_{10} \cup \operatorname{supp}u_{20}=\Omega,\quad \operatorname{supp}u_{10}\cap\operatorname{supp}u_{20}=\omega_0,
\end{equation}
where $\omega_0\subset\Omega$ is the contact hypersurface, one or several points in our one-dimensional simulations. The theory predicts that the initial segregation is kept for later times, that is, no mixing between species is possible.

Under our modeling point of view, the initial data segregation can be interpreted
as if the process of species splitting (tumor cells and normal cells) took place in coincidence with that of niches segregation: tumor tissue on one side of $\omega_0$, and normal tissue on the other side.

For the one-dimensional model, the results of Bertsch et al. \cite{bertsch2010}, or the results proven in the present article include this special case of initial data. 

For the multi-dimensional model, existence of solutions of problem (D) with initial data satisfying \fer{ci.id} is proven in \cite{bertsch2012,gsv}.  These multi-dimensional results are based on the use of suitable Lagrangian formulations of the problem but, as suggested by its parabolic-hyperbolic nature evidenced in the auxiliary formulation introduced in \cite{bertsch2010}, see Problem (P)$_B$, this construction lacks of uniqueness, as proven in  \cite{gsv}. 

Thus, the viscosity approximation used to prove the existence of solutions in the one-dimensional case \cite{Galiano2014} could be a method to select one of the infinitely many solutions obtained by the Lagrangian approach. This, however, is just a conjecture.

Finally, notice that the case of general initial data in the multi-dimensional framework  is still an open problem.

In this example, we consider a particular situation of the contact-inhibition problem where a 
%$BV(\O)$ 
space  discontinuous  explicit solution of problem~(ND) may be computed in terms of the
Barenblatt explicit solution of the porous medium equation, the Heaviside function and the
trajectory of the contact-inhibition point. To be precise, we construct a solution to the
problem, 
\begin{align}
 %\label{prob:pert.ec}
& \partial_t u_{i} - \p_x(u_i \p_x(u_1+u_2)) =0&& \text{in }Q_T,&\label{eq:s1}\\
&u_i \p_x(u_1+u_2) = 0 &&\text{on } \Gamma_T,& \label{eq:s2}
\end{align}
with, for $x,x_0\in\O=(-L,L)$,
\begin{equation}
\label{def:us}
 u_{10}(x)=H(x-x_0)B(0,x), \quad u_{20}(x)=H(x_0-x)B(0,x).
\end{equation}
Here, $H$ is the Heaviside function and $B$ is the Barenblatt solution of the porous medium
equation corresponding to the initial datum $B(-t^*,\cdot)=\delta_0$ (Dirac delta function), i.e.
\begin{equation*}
%\notag \label{barenblatt}
 B(t,x)= 2 (t+t^*)^ {-1/3} \big[1-\frac{1}{12}x^ 2(t+t^*)^ {-2/3}\big]_+,
\end{equation*}
with $t^* >0$.
For simplicity, we consider problem \eqref{eq:s1}-\eqref{def:us} for $T>0$ such that
$\rho(T)<L^ 2$, with $\rho(t)=\sqrt{12}(t+t^*)^ {1/3}$, so that $B(t,\pm L)=0$ for all $t\in[0,T]$.
The point $x_0$ is the initial contact inhibition point, for which we assume $|x_0|<\rho(0)$,
i.e. it belongs to the interior of the support of $B(0,\cdot)$, implying that the initial
mass of both populations is positive.

\begin{theorem}\label{th.easy}
 The $L^\infty(0,T;BV(\O))$ functions
\begin{equation}
\label{def.uib}
 u_1(t,x)=H(x-\eta(t))B(t,x),\quad u_2(t,x)=H(\eta(t)-x)B(t,x),
\end{equation}
with $\eta(t)=x_0(1+ \frac{t}{t^*})^{1/3}$, are a weak solution of problem \eqref{eq:s1}-\eqref{def:us}
in the following sense: 
 For any $\vfi \in C^1(0,T;H^1(\O))$ with $\vfi(T,\cdot)=0$ in $\O$, 
 \begin{align*}
 \int_{Q_T} u_i\p_t \vfi = \int_{Q_T} (u_i\p_x u ) \p_x\vfi- \int_{-L}^{L} u_{i0}\vfi(0,\cdot).
 %\label{def.wd}
\end{align*}
\end{theorem}

\begin{figure}[h] 
\centering 
{\includegraphics[width=3cm,height=3cm]{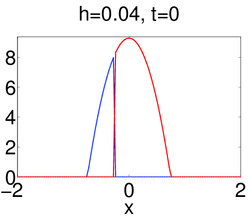}}
{\includegraphics[width=3cm,height=3cm]{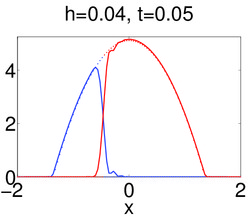}}
{\includegraphics[width=3cm,height=3cm]{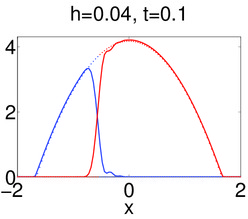}}
{\includegraphics[width=3cm,height=3cm]{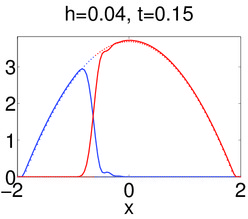}}\\
{\includegraphics[width=3cm,height=3cm]{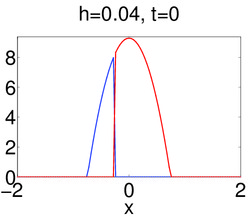}}
{\includegraphics[width=3cm,height=3cm]{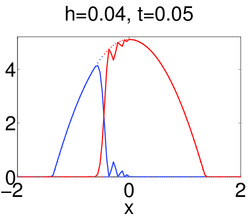}}
{\includegraphics[width=3cm,height=3cm]{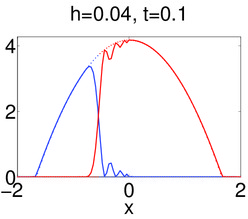}}
{\includegraphics[width=3cm,height=3cm]{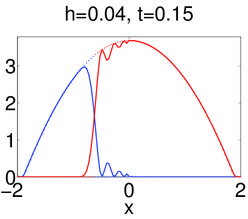}}\\
{\includegraphics[width=3cm,height=3cm]{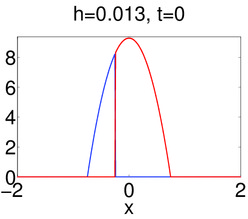}}
{\includegraphics[width=3cm,height=3cm]{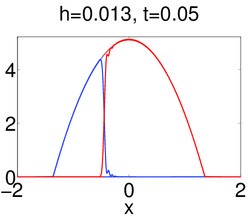}}
{\includegraphics[width=3cm,height=3cm]{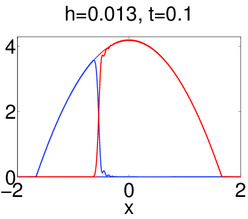}}
{\includegraphics[width=3cm,height=3cm]{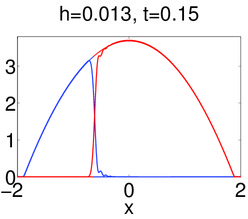}}\\
{\includegraphics[width=3cm,height=3cm]{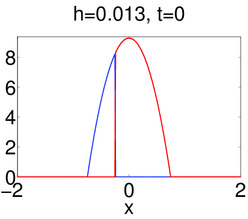}}
{\includegraphics[width=3cm,height=3cm]{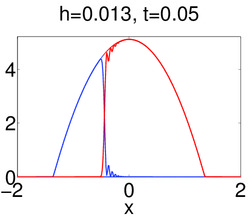}}
{\includegraphics[width=3cm,height=3cm]{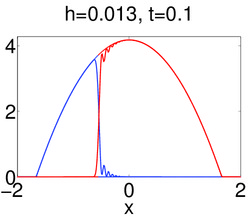}}
{\includegraphics[width=3cm,height=3cm]{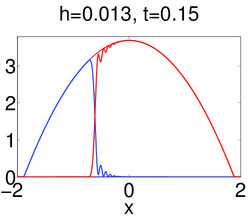}}\\
{\includegraphics[width=3cm,height=3cm]{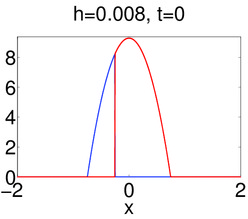}}
{\includegraphics[width=3cm,height=3cm]{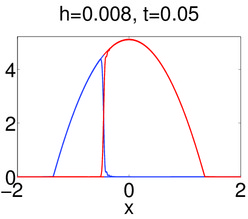}}
{\includegraphics[width=3cm,height=3cm]{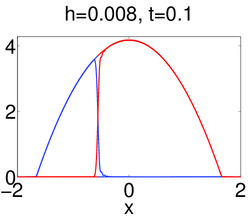}}
{\includegraphics[width=3cm,height=3cm]{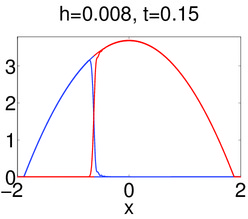}}\\
{\includegraphics[width=3cm,height=3cm]{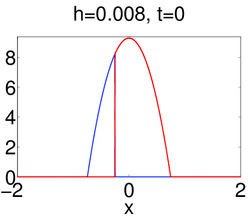}}
{\includegraphics[width=3cm,height=3cm]{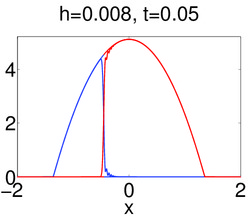}}
{\includegraphics[width=3cm,height=3cm]{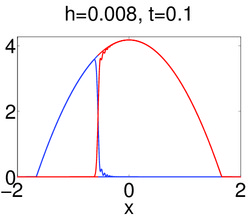}}
{\includegraphics[width=3cm,height=3cm]{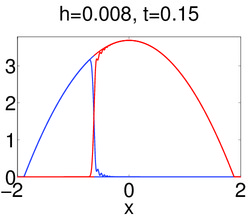}}\\
\caption{  Experiment 2. Exact solution (dots) against approximated solutions (solid lines). 
Odd rows correspond to the approximated solution of problem (P)$_\delta$ for several time slices, $t$, while even rows correspond to the solution of problem (P)$_B$.
Rows correspond to different mesh sizes, captured by parameter $h$. Mind the different vertical scales among time slices.}
\label{fig_barenblatt_sol} 
\end{figure}

With this experiment, we check whether the condition $u_1+u_2>0$ in $Q_T$ needed for proving the convergence of the solutions of problems (P)$_\delta$ and (P)$_B$ to a solution of problem (D) is necessary or not. In particular, we consider $L=2$,  $x_0=-0,25$ and $t^*=0.01$, so that $u_1(t,\cdot)+u_2(t,\cdot)=0$ in some regions of $\O=(-2,2)$ for $t<T=0.15$, see Figure~\ref{fig_barenblatt_sol}.
%
%Thus, we consider the following data: $\O=(-2,2)$, $x_0=-0,25$, $t^*=0.01$, so that $u_1(t,\cdot)+u_2(t,\cdot)=0$ in some regions of $\O$ for $t<T=0.15$, see Figure~\ref{fig_barenblatt_sol}.
%
%We performed simulations for getting numerical approximations of problems (P)$_\delta$ and (P)$_B$
%with regularization parameters $\delta=h^2$, $\delta_B=2h^2$, respectively, where $h$ is the mesh size.
%For $\lambda_\eps$ and $\Lambda_\eps$, mainly avoiding negative values of the discrete solutions, we set $\eps=1.e-10$.  The tolerance for the fixed point stopping criterion is fixed as $\text{tol}=1.e-8$.
%Finally, the time step was empirically chosen as $\tau=1.e-4$. In all the examples, we observed good convergence properties of the fixed point algorithm, reaching the prescribed tolerance in less than ten iterations. 
%
In fact, we already have initial conditions such that  $u_{10}+u_{20}=0$ in some regions of $\O$. According to this, and on the contrary of problem (P)$_\delta$, the initialization of the data for problem (P)$_B$ requires special care since by construction we have $r_0=u_{10}/(u_{10}+u_{20})$ in $\O$. In this sense, among several choices to define $r_0$ % when $u_{10}+u_{20}$ vanishes, 
the following gave us the best results: Considering the perturbation, for $\gamma>0$,
\begin{equation*}
u_{10}^\gamma =\left\{\begin{array}{ll}
                \gamma+u_{10}  & \text{if }x\leq x_0, \\
                0 & \text{if }x > x_0 ,
               \end{array}
               \right.\quad 
u_{20}^\gamma =\left\{\begin{array}{ll}
                0 & \text{if }x < x_0 ,\\
                \gamma+u_{20}  & \text{if }x\geq x_0 ,
               \end{array}
               \right.\               
\end{equation*}
we get $u_{10}^\gamma+u_{20}^\gamma =\gamma +u_{10}+u_{20} >0$. Then, $r_{0}^\gamma(x)=1$ if $x<x_0$ 
and $r_{0}^\gamma(x)=0$ if $x>x_0$, and thus, taking $\gamma\to 0$ we obtain the same definition  for $r_0$.

\begin{figure}[h] 
\centering 
{\includegraphics[width=4cm,height=4cm]{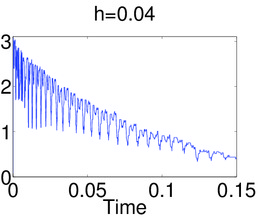}}
{\includegraphics[width=4cm,height=4cm]{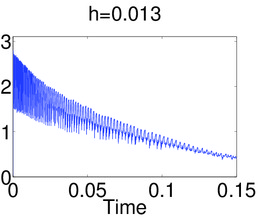}}
{\includegraphics[width=4cm,height=4cm]{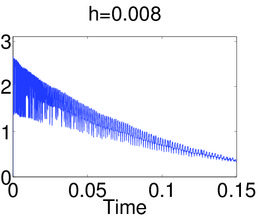}}\\
{\includegraphics[width=4cm,height=4cm]{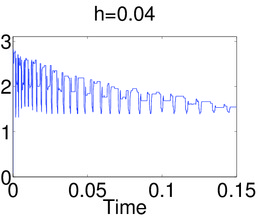}}
{\includegraphics[width=4cm,height=4cm]{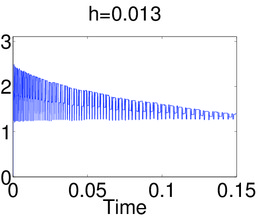}}
{\includegraphics[width=4cm,height=4cm]{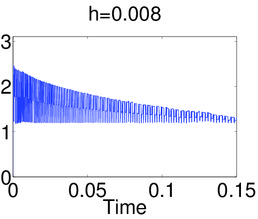}}\\
{\includegraphics[width=4cm,height=4cm]{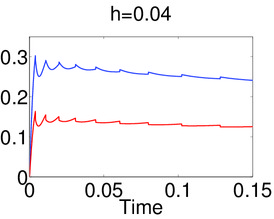}}
{\includegraphics[width=4cm,height=4cm]{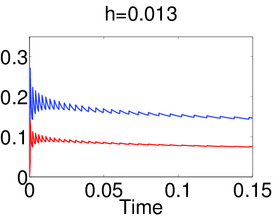}}
{\includegraphics[width=4cm,height=4cm]{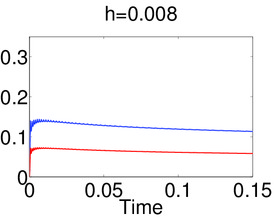}}\\
{\includegraphics[width=4cm,height=4cm]{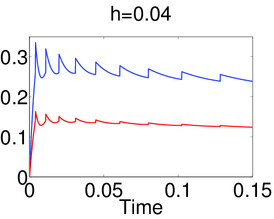}}
{\includegraphics[width=4cm,height=4cm]{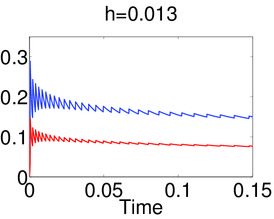}}
{\includegraphics[width=4cm,height=4cm]{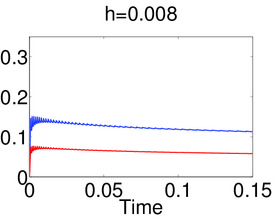}}
\caption{ Experiment 2. Two first rows correspond to the oscillation measure $\textrm{osc}(u)(t)$ of the solutions approximated by (P)$_\delta$ and (P)$_B$, respectively. 
Two last rows correspond to the relative error between the exact solution and the solution approximated by (P)$_\delta$ and (P)$_B$, respectively. 
Columns correspond to different mesh sizes, captured by parameter $h$.}
\label{fig_barenblatt_error} 
\end{figure}

In Figures~\ref{fig_barenblatt_sol} and \ref{fig_barenblatt_error} we plot the discrete solutions 
obtained from both approximation schemes together with the explicit solution, and the corresponding oscillation measure.
Just as happened in the previous example, in all the cases the solutions corresponding to problem (P)$_B$ produce more oscillations than those of problem (P)$_\delta$, and solutions of  (P)$_\delta$ attenuate the oscillations when $t$ increases whereas those corresponding to solutions of  (P)$_B$ are always above some threshold value. Moreover, in this example the oscillations come up around the contact inhibition point. 

It is also interesting to notice that, since this example has an explicitly known solution, we can check the convergence of the numerical solutions to the exact one. Indeed, the relative $L^2$ errors of both solutions are similar, and decreasing with respect to the mesh size, see Figure~\ref{fig_barenblatt_error}, last two rows.

\section{Proofs of the theorems}\label{sec:proof}

\subsection{Proof of Theorem~\ref{th.1}}

We divide the proof in several steps.

\no\emph{Step 1. An approximation problem. }We consider problem \Pd~
for some non-negative $\ud_{i0}\in C^1(\bar\O)$ with $\ud_{i0}\to u_{i0}$ strongly in $BV(\O)$, $ \min(u_0)\leq \ud_{10}+\ud_{20}\leq \max(u_0)$ in $\O$, and 
 \begin{align}
  \label{ui0.b}
  \sqrt{\delta}\nor{\p_x \ud_{i0}}_{L^2(\O)} \text{  uniformly bounded.}
 \end{align}

This nonlinear viscosity regularization allows us to use the results in \cite{ggj2,chen} to deduce the existence of a  solution $(\ud_1,\ud_2)$, with non-negative $\ud_i \in L^2(0,T;H^1(\O))\cap W^{1,4/3}(0,T;(W^{1,4}(\O))')$,  satisfying \Pd~ in the weak sense:
\begin{align}
 \int_0^T < \p_t \ud_i,\vfi> & +\int_{Q_T} \big(\ud_i\p_x \ud + \ud_iq +\frac{\delta}{2} \p_x(\ud_i \ud) \big) \p_x\vfi \nonumber \\
 & = \int_{Q_T} f_i(\ud_1,\ud_2)\vfi,
 \label{def.wd}
\end{align}
for all $\vfi\in L^{4}(0,T;W^{1,4}(\O))\cap L^\infty(Q_T)$. Here, $<\cdot,\cdot>$ denotes the duality product 
$(W^{1,4}(\O))'\times W^{1,4}(\O)$.

Adding the corresponding equations of problem \Pd, for $i=1,2$, we find that $\ud=\ud_1+\ud_2$ satisfies, in a weak sense, the following problem \Ps
\begin{align*}
 &\p_t \ud  -(1+\delta) \p_x (\ud \p_x \ud) + \p_x(q\ud )
 = F(\ud_1,\ud_2)& & \text{in } Q_T, & \\
  &\p_x \ud  =0 & & \text{on } \Gamma_T, & \\
  &\ud(0,\cdot)=u_{10}^\delta+u_{20}^\delta&& \text{in } \O, & 
\end{align*}
where $F(s_1,s_2)= f_1(s_1,s_2)+f_2(s_1,s_2)$. Using Theorem~2.2 (step 5) of \cite{dgj2001} and the 
pointwise bounds for $\ud_{i0}$, we get 
\begin{align}
\label{aux1}
\min(u_0)e^{-\lambda t}\leq \ud (t,\cdot) \leq \max(u_0)e^{\lambda t}\qtext{in }\O,\qtext{for a.e. }t\in (0,T),
\end{align}
for some $\lambda>0$ independent of $\delta$. 
Thus, since $\ud_i\geq 0$ in $\O$, we also get
\begin{align}
\label{est.uiinf}
0\leq \ud_i (t,\cdot) \leq \max(u_0)e^{\lambda t}\qtext{in }\O,\qtext{for a.e. }t\in (0,T).
 \end{align}
The uniform estimate \fer{aux1} together with the regularity of $u_0$ and $q$, and the compatibility conditions satisfied by these functions, allow us to apply the general theory for evolution quasilinear equations to deduce that the unique solution of \Ps~is Lipschitz continuous in $Q_T$, and that, for any $1\leq p \leq\infty$,
\begin{align}
\label{reg.u}
 \nor{\p_t \ud}_{L^p(Q_T)},~\nor{\p_{xx}\ud}_{L^p(Q_T)}  \text{ are  uniformly bounded,}
\end{align}
see \cite[Theorem 7.20]{lieberman}.

Returning to problem \Pd, we may see that the additional regularity $\ud_i\in L^\infty(Q_T)$ given by \fer{est.uiinf} allows us to expand the sense of weak solution to the regularity 
\begin{align*}
 \ud_i \in L^2(0,T;H^1(\O))\cap H^1(0,T;(H^1(\O))') \cap L^\infty(Q_T),
\end{align*}
and for test functions $\vfi\in L^{2}(0,T;H^1(\O))$. Using the arguments of the proof of \cite[Lemma 3]{Galiano2014}, we easily get 
\begin{align}
\label{est.ui0}
 \sqrt{\delta}\nor{\p_x \ud_i}_{L^2(Q_T)}, ~\nor{\p_t \ud_i}_{L^2(0,T;(H^1(\O))')} \qtext{are uniformly bounded.}
\end{align}
Moreover, the above mentioned regularity  $\p_{xx}\ud\in L^\infty(Q_T)$ allows us to express  \fer{def.wd} as, for all $\vfi\in L^{2}(0,T;H^1(\O))\cap L^\infty(Q_T)$, 
\begin{align*}
  \int_0^T  <\p_t \ud_i , \vfi > & + \frac{\delta}{2} \int_{Q_T} \ud \p_x\ud_i \p_x \vfi = 
   \int_{Q_T} \psi\vfi,
\end{align*}
with 
\begin{align*}
 \psi = (1+\frac{\delta}{2}) ( \p_x \ud_i \p_x \ud + \ud_i \p_{xx}\ud ) 
  + \ud_i \p_x q +\p_x \ud_i q  + f_i(\ud_1,\ud_2) .
\end{align*}
Since $\psi\in L^2(Q_T)$, we have that $\ud_i$ satisfies 
\begin{align*}
 \p_t \ud_i -\frac{\delta}{2} \p_{x}(\ud \p_x\ud_i) =\psi \in L^2(Q_T),
\end{align*}
with $\ud\in C^{0,1}(\bar Q_T)$ bounded away from zero, see \fer{aux1}. Thus, we have 
$\p_t \ud_i,~\p_{xx} \ud_i \in L^2(Q_T)$, implying $\p_x \ud_i\in L^4(Q_T)$. We then deduce that, in fact, $\psi \in L^4(Q_T)$, improving in this way the regularity of $\ud_i$. A boot-strap argument allows us to deduce 
\begin{align}
 \label{reg.ui1}
 \p_t \ud_i,~\p_{xx} \ud_i \in L^\infty(Q_T), 
\end{align}
implying the Lipschitz continuity of $\ud_i$ (with norm depending on $\delta$).
\bigskip

\no\emph{Step 2. Uniform estimates for $\rd=\ud_1/\ud$.}  For clarity, we omit the super-index $\delta$ in this Step. Due to the regularity \fer{reg.ui1}, the derivatives $\p_t r$, 
$\p_x r$, and $\p_{xx} r$ are well defined as $L^\infty(Q_T)$ functions. After some computations, we obtain that $r$ satisfies
\begin{align}
 &\p_t r -\frac{\delta}{2}u \p_{xx} r - ( q+(1+2\delta) \p_x u) \p_x r  = G(u,r)
 & & \text{in } Q_T, &  \label{eq.r} \\
  &\p_x r  =0 & & \text{on } \Gamma_T, & \nonumber \\
  &r(0,\cdot)=u_{10}/u_0&& \text{in } \O,\nonumber & 
\end{align}
with, 
\begin{align*}
 G(u,r)=   \frac{1-r}{u} f_1(ur,u(1-r))- \frac{r}{u}f_2(ur,u(1-r)).
\end{align*}
Since $\partial_u G$ and $\partial_r G$ are bounded, and $u$ and $r$ are Lipschitz continuous, we deduce $\p_x G(u,r)\in L^\infty(Q_T)$.  Therefore, the solution $r$ of \fer{eq.r} is  regular enough to allow us to differentiate equation \fer{eq.r} with respect to $x$, for a.e. $(t,x)\in Q_T$. We then multiply the differentiated equation by $\xi= \p_x r /\sqrt{\delta+ (\p_x r)^2}$ and integrate by parts to get, term by term,
\begin{align*}
 & \int_{Q_T} \p_x \p_t r \xi  = \int_\O \sqrt{\delta+ (\p_x r)^2} \Big|_{0}^{T} ,\\
   -\frac{\delta}{2} & \int_{Q_T} \p_x\big(u\p_{xx}r \big)\xi  = 
 \frac{\delta^2}{2} \int_{Q_T} \frac{u(\p_{xx} r)^2}{(\delta+ (\p_x r)^2)^{3/2}} ,\\
  -& \int_{Q_T} \p_x \big(  ( q+(1+2\delta) \p_x u) \p_x r \big) \xi = 
\delta \int_{Q_T}\frac{\p_x ( q+(1+2\delta) \p_x u)}{\sqrt{\delta+ (\p_x r)^2}} ,\\
& \int_{Q_T} \p_x G(u,r) \xi = \int_{Q_T} (\p_u G(u,r) \p_x u +\p_r G(u,r) \p_x r ) \xi.
\end{align*}
Therefore, thanks to the regularity $\p_x q \in L^1(Q_T)$ and to the uniform bounds for $\p_x u$ and $\p_{xx} u$ in $L^1(Q_T)$ we obtain
\begin{align*}
 \int_\O \sqrt{\delta+ (\p_x r)^2} \Big|_{0}^{T} \leq  
 c_1 + c_2 \int_{Q_T} \abs{\p_x r},
\end{align*}
with $c_1,c_2$ independent of $\delta$. We then deduce from Gronwall's lemma that 
\begin{align}
\label{est.rx1}
 \sup_{[0,T]}\int_{\O} \abs{\p_x r} \qtext{is uniformly bounded with respect to }\delta.
\end{align}

We finish this step by showing a uniform bound for $\p_t r$. From equation \fer{eq.r} we get
 \begin{align}
 \label{est.rt}
  \int_{Q_T} \abs{\p_t r} \leq \frac{\delta}{2}\nor{u}_{L^\infty(Q_T)} \int_{Q_T}\abs{\p_{xx} r}
  + c_3 \int_{Q_T} \abs{\p_x r }+  \int_{Q_T} \abs{G(u,r)},
 \end{align}
with $c_3=\nor{q}_{L^\infty(Q_T)} +(1+2\delta) \nor{\p_x u}_{L^\infty(Q_T)}$. We have already shown that the two last terms of the right hand side of \fer{est.rt} are uniformly bounded. For estimating the first term, we multiply equation \fer{eq.r} by $\delta \p_{xx}r$ and integrate by parts. We obtain
\begin{align*}
 \min(u) \frac{\delta^2}{2}\int_{Q_T} \abs{\p_{xx} r}^2 \leq & \frac{\delta}{2}\int_\O \abs{\p_x r(0,\cdot)}^2 +  \frac{\delta(1+2\delta)}{2} \int_{Q_T} (\abs{\p_{xx} u}+\abs{\p_{x} q})  \abs{\p_{x} r}^2 \\
 &+ \delta 
 \int_{Q_T} \abs{G(u,r)}\abs{\p_{xx} r}.
\end{align*}
Using property \fer{ui0.b}, we find that the first term of the right hand side is uniformly bounded. The third term may be handled by H\"older's and Young's inequality to get 
\begin{align*}
 \min(u) \frac{\delta^2}{4}\int_{Q_T} \abs{\p_{xx} r}^2 \leq & c_4 + \frac{\delta(1+2\delta)}{2} \int_{Q_T} (\abs{\p_{xx} u}+\abs{\p_{x} q})  \abs{\p_{x} r}^2 ,
\end{align*}
with $c_4$ independent of $\delta$. Finally, using the uniform lower bound for $u$ in $Q_T$, the uniform bound for $\p_{xx} u$ in $L^2(Q_T)$, the regularity $\p_x q \in L^2(Q_T)$ and applying the last argument of point (v) of the proof of Lemma 3.1 of \cite{bertsch2010}, we deduce 
that $\delta \nor{\p_{xx} r}_{L^1(Q_T)}$ is uniformly bounded,
% \begin{align*}
%  \delta \nor{\p_{xx} \rd}_{L^1(Q_T)} \qtext{is uniformly bounded,}
% \end{align*}
and therefore, from \fer{est.rt} we also deduce 
\begin{align}
\label{est.rx2}
  \nor{\p_{t} r}_{L^1(Q_T)} \qtext{is uniformly bounded with respect to }\delta.
\end{align}

\no\emph{Step 3. Passing to the limit $\delta\to0$. } 
From the uniform estimates \fer{est.uiinf} and \fer{est.ui0},  we deduce the existence of 
$u_i\in  L^\infty(Q_T) \cap L^2(0,T;(H^1(\O))') $ such that (for a subsequence, not relabeled)
\begin{align*}
 & \ud_i\wto u_i \qtext{weakly * in } L^\infty(Q_T),\\
 & \p_t \ud_i\wto \p_t u_i \qtext{weakly in } L^2(0,T;(H^1(\O))'),\\
 & \delta \p_x \ud_i\to 0 \qtext{strongly in } L^2(Q_T),
\end{align*}
and from \fer{reg.u}, we also deduce the existence of $u\in C^{0,1}(\bar Q_T)\cap L^2(0,T;H^2(\O))$ such that 
\begin{align*}
  & \ud\to u \qtext{uniformly in }  C^{0,1}(\bar Q_T),\\
 & \p_x \ud\to \p_x u \qtext{strongly in } L^2(Q_T).
\end{align*}
In particular, the weak and strong convergences of $\ud_i$ and $\ud$ imply $u=u_1+u_2$.
Using the uniform estimates \fer{est.rx1} and \fer{est.rx2} we deduce the existence of $r\in L^\infty(0,T;BV(\O))\cap BV(0,T;L^1(\O))$ such that, for $1\leq p<\infty$,
\begin{align*}
  & \rd\to r \text{ strongly in } L^p(Q_T). 
\end{align*}
Since $\rd=\ud_1/\ud$ the weak and strong convergences of $\ud_i$ and $\ud$ in $L^2(Q_T)$, respectively, imply $r=u_1/u$. Then, the strong convergences of $\rd$ and $\ud$ in $L^p(Q_T)$ and $L^\infty(Q_T)$ imply
\begin{align*}
 \ud_1=\rd\ud \to ru=u_1 \text{ strongly in } L^p(Q_T). 
\end{align*}
Similarly, we obtain 
\begin{align*}
 \ud_2=(1-\rd)\ud \to (1-r)u=u_2 \text{ strongly in } L^p(Q_T).
\end{align*}
Finally, using again the uniform estimates \fer{est.rx1} and \fer{est.rx2}, we also obtain 
\begin{align*}
 \nor{\p_x \ud_i}_{L^\infty(0,T;L^1(\O))},\nor{\p_t \ud_i}_{L^1(Q_T)}\text{ uniformly bounded,}
\end{align*}
from where we, additionally, deduce 
\begin{align*}
 u_i\in L^\infty(0,T;BV(\O))\cap BV(0,T;L^1(\O)).
\end{align*}
Thus, the passing to the limit $\delta\to 0$ in \fer{def.wd} is justified, and so we identify the limit $(u_1,u_2)$ as a weak solution of problem (D). $\Box$

\subsection{Proof of Theorem~\ref{th.easy}}

\noindent\emph{Proof. }
Let $H_\epsilon$ be the regularization of the Heaviside function taking the values
$\left\{0,\frac{1}{2}(1+x/\epsilon),1\right\}$ in the intervals $(-L,-\epsilon)$,
 $(-\epsilon,\epsilon)$ and $(\epsilon,L)$, respectively, for $\epsilon>0$ small.

 Define the functions $u_i^\epsilon$ as in \fer{def.uib}, with $H$ replaced by $H_\epsilon$, and 
 let $\varphi\in H^ 1(Q_T)$, with $\vfi(T,\cdot)=0$ in $\O$.  Using $\varphi_\epsilon
(t,x)=\varphi(t,x) H_\epsilon(x-\eta(t)) $ as a test function in the weak formulation of
\fer{eq:s1}-\fer{eq:s2} corresponding to the initial datum $B(0,\cdot)$, that is, 
the problem satisfied by the Barenblatt solution in $Q_T$, we obtain
\begin{align}
\label{weak.B}
   \int_{Q_T}  H_\epsilon(x-\eta(t)) B(t,x) (\p_t\varphi(t,x) & - \p_x B(t,x) \p_x\varphi(t,x) )dxdt \\
 & +\int_{-L}^ {L}  H_\epsilon(x-x_0) B(0,x) \varphi(0,x)dx 
    =I^ 1_\epsilon, \nonumber
\end{align}
with
\begin{eqnarray*}
I^ 1_\epsilon
  = \int_{Q_T} \varphi(t,x) B(t,x) H'_\epsilon(x-\eta(t)) \big(\eta'(t) + \p_x B(t,x)\big) dx dt.
\end{eqnarray*}
Since $|x_0|<\rho(0)$, we have $\abs{\eta(t)}<L^2-\epsilon$, for $\epsilon$ small enough and
$t\in(0,T)$. Using the explicit expression of $\p_x B$ and $\eta'$ we deduce
\begin{equation*}
%\notag \label{eq.I} 
I_\epsilon^1=-\frac{1}{6\epsilon}\int_{0}^T\int_{-\epsilon}^\epsilon y
\varphi(t,y+\eta(t))B(t,y+\eta(t))dydt.
\end{equation*}
Since $\varphi$ and $B$ are uniformly bounded in $L^\infty(Q_T)$, we obtain

\begin{equation}
\label{est.I}
 |I^ 1_\epsilon|\leq C \epsilon ,
\end{equation}
with $C>0$ independent of $\epsilon$. The computation using
$\varphi(t,x)H_\epsilon(\eta(t)-x)$ as test function gives similar results for some
$I_\epsilon^ 2$ satisfying the same estimate \eqref{est.I} than $I_\epsilon^ 1$. 
Since, by definition, $u_1^\epsilon+u_2^\epsilon=B$, we obtain
\begin{equation*}
\left|\int_{Q_T}  u_{i}^ \epsilon (\p_t \varphi -  \p_x(u_1^ \epsilon+u_2^
\epsilon ) \p_x\varphi) + \int_{-L}^ {L}  (u_i^{\epsilon}
\varphi)(0,\cdot) \right| \leq C\epsilon
\end{equation*}
for all $\varphi\in H^ 1(Q_T)$.
Thus, using that  $u_i^ \epsilon \to u_i$ uniformly 
in $L^ \infty (Q_T)$, and performing the limit $\epsilon\to 0$ in \fer{weak.B}, we deduce that $(u_1,u_2)$ is a solution of problem \fer{eq:s1}-\fer{def:us}.
$\Box$
 
\begin{remark}
It is not difficult to extend the above construction to problem~(ND) for non-vanishing $q$ and $f_i$. Indeed, reconsider the solution $u$ of problem \fer{p.u}
 and the corresponding approximations $u_i^\epsilon$ given in the proof of the previous theorem. Then, to handle the integrals $I_\epsilon^i$, we first observe that
for $\epsilon\to0$ we get
\begin{eqnarray*}
 I^ 1_\epsilon \to -\int_0^T \varphi(t,\eta(t)) u(t,\eta(t)) \big(\eta'(t) - b(t,\eta(t))\big) dt,
\end{eqnarray*}
with $b = -(\p_x u + q).$ Therefore, if the ODE problem
\begin{equation}
\label{eq.eta2}
\left\{
 \begin{array}{ll}
  \eta'(t)=b(t,\eta(t)) & \text{for }t\in (0,T), \\
  \eta(0)=x_0,
 \end{array}
\right.
 \end{equation}
is solvable, a solution for problem (ND) may be constructed as before. Typical conditions on
$b$ for \eqref{eq.eta2} to be solvable (in the one-dimensional case) are given in terms of 
spatial Lipschitz continuity. That is, the solution of problem \fer{p.u} must satisfy $\p_{xx}u\in L^\infty(Q_T)$, which is true for smooth data.
\end{remark}

\section{Conclusions}\label{sec:conclusions}

The generalization of the dynamical system model for the splitting - differentiation process of populations \cite{sanchez-palencia} to the space dependent situation leads to a family of cross-diffusion PDE problems including several well known segregation models.

In a previous work \cite{Galiano2012}, we analyzed the case in which the cross-diffusion is of the type introduced by Shigesada et al. \cite{SKT79}, a model well studied in the literature, see \cite{Galiano2014} and its references. In this article, we turned to another important case, the model of Busenberg and Travis \cite{busenberg83}, in which the segregation effects are stronger, as motivated in \cite{Galiano2014}. 

The model of Busenberg and Travis was studied by Bertsch et al. in a series of papers \cite{bertsch85,bertsch2010,bertsch2012}   
focusing in a special case of initial data arising from Tumor modeling which gives rise to the so-called contact-inhibition problem.   

From the analytical side, we have produced a proof of existence of solutions conceptually simpler than previous proofs \cite{bertsch2010,gsv}, since ours derives from a direct parabolic regularization of the problem meanwhile previous involved a change of unknowns rendering the problem to a parabolic-hyperbolic formulation. 

Although difficult to ensure, we believe that our approach also gives  
a way to select a unique \emph{natural} solution of the problem as a limit of vanishing viscosity solutions of the regularized parabolic problems, tackling thus the non-uniqueness issue of the 
parabolic-hyperbolic formulation \cite{gsv}. 

From the numerical side, we have introduced a Finite Element discretization of both formulations 
and compared the results in a series of experiments (only two of them showed in the article). 
The general observation is that our approach is always more stable in the tricky regions where the solutions exhibit discontinuities. 

Finally, there are  evidences in our numerical experiments showing that the main data 
restriction of the existence theorems, the positivity of the initial total mass, is not 
a necessary condition. Thus, future research will focus in the replacement of this condition, as
well as in the generalization to the multi-dimensional case.

%\bibliographystyle{elsarticle-num} 

%\bibliographystyle{ecology}
%\bibliographystyle{spbasic}
%\bibliography{discrete-to-continuous}

\end{document}